\input amstex
\documentstyle{amsppt}

\refstyle {A}
\magnification=\magstep1 

\hoffset .25 true in
\voffset .2 true in

\hsize=6.1 true in
\vsize=8.5 true in

\def\op{\operatorname}
\def\s{\sigma}

\def\D{\Cal D}
\def\C{\Bbb C}
\def\E{\Cal E}
\def\g{\frak g}
\def\ot{\otimes}
\def\X{\bold X}

\def\O{\Omega}
\topmatter
\title
On quantum Jacobi identity
\endtitle

\author Maxim Vybornov \endauthor

\address Department of Mathematics, Yale University, 
New Haven, CT 06520
\endaddress

\email mv\@math.yale.edu \endemail

\date November, 1996\enddate

\abstract
In this paper we present our variant of quantum antisymmetry
and quantum Jacobi identity. 
\endabstract
\endtopmatter

\document

\subhead  1. Introduction \endsubhead   The notion of 
``quantum Lie algebra,'' related to the notion of quantum group,
has recently been investigated by  many researchers 
(see Bibliography on quantum Lie algebras.)
It seems, however,  that many important problems are still open.

Our goal in the subject is to reconstruct the category $\E$
(of finitely generated representaions of quantum group
corresponding to complex simple Lie algebra $\g$) from
such data as the quantization of the
adjoint representation $V^h$, the quantum Casimir tensor $t^h$,
the quantum Lie bracket $T^h$, the quantum commutativity
morphism $\check R^h$, etc i.e. we would like 
to present the category $\E$ as the category
of modules over quantum Lie algebra.

As a first step towards the goal  we present
in this paper our version of what may be called
``quantum antisymmetry'' and ``quantum Jacobi identity.''
It turns out that the infinitesimal 
deformation resembles the classical Jacobi identity.
\medskip

\subhead  2. Notation \endsubhead  
In this paper $1$ stands for the identity morphism of some
object. Each time it is either clear  which object we mean,
or we specify the object.   

Let $\g$ be a complex 
semisimple Lie algebra, and 
$t\in \g\otimes \g$ be the invariant symmetric tensor
corresponding to the Killing form on $\g$.
Following  notations and terminology of Kazhdan and 
Lusztig ({\bf [6]}),
we take   $\D$ and $\E$ to  be Drinfeld's categories associated to
$\g$. We notice that the category $\D$ has 
two structures of braided 
tensor category: a trivial structure $(\D, \ot, \s, 1\ot 1\ot 1)$,
where $\s:U\ot W @>>> U\ot W$ is just the simple 
transposition of factors in $U\ot W$, and a non-trivial structure 
$(\D, \ot, \check R=\s\circ R$ \linebreak 
$=\s\circ e^{\sqrt{-1}\pi ht},
\Phi= \Phi_{KZ}(t_{12}, t_{23}))$.
Let $\X:\D @>>>\E$ be the equivalence functor of 
Kazhdan and Lusztig.
For two objects $U$ and $W$ in $\D$, let
$M_{U,W}: \X(U)\ot \X(W) @>>> \X(U\ot W)$ be the ``natural 
transformation'' morphism providing the tensor structure of 
the functor $\X$ (see 25.7. of  {\bf [6]}.)
Transforming the trivial structure in $\D$ via the functor $\X$
we will get the second structure of the braided tensor category
in $\E$. We will denote the corresponding commutativity morphism
by $\s^h$ and the corresponding associativity morphism by $\Psi$.
We will study the category $(\E, \ot, \s^h, \Psi)$ in more detail. 

In the category $\D$ we have a distinguished 
object $V=\g [[h]]$ and a distinguished morphism 
$[\ ,\ ]=T:V\otimes V @>>> V$, corresponding to the Lie bracket
on $\g$, and a distinguished morphism $\check R:V\ot V @>>> V\ot V$,
$\check R=\s\circ R_{V,V}$, where $R_{V,V}$ is 
the representation in 
$V\ot V$
of the universal $R$-matrix $R=e^{\sqrt{-1}\pi ht}$.
The morphisms $t^h$, $T^h$, and $\check R^h$ in $\E$ are 
defined as
$$
t^h=M_{V,V}^{-1}\circ\X(t)\circ M_{\C[[h]],\C[[h]]},
\qquad
T^h=\X(T)\circ M_{V,V},
\qquad
\check R^h=M_{V,V}^{-1}\circ \X(\check R)\circ M_{V,V}
$$
where $t$ is considered as a morphism 
$t: \C[[h]]\ot\C[[h]] @>>> V\ot V$.
We denote by $V^h=\X(V)$ and we set
$$
\matrix
\format\l & \l\\
\O=-(T\ot T)\circ (1_V\ot t\ot 1_V) & \\
\O^h=M_{V,V}^{-1}\circ \X(\O)\circ M_{V,V}  & \\
\s^h=M_{V,V}^{-1}\circ \X(\s)\circ M_{V,V}
=\check R^h\circ e^{-\sqrt{-1}\pi h 
\O^h}
\endmatrix
$$
where $\O$ is considered as a morphism
$V\ot V @>>> V\ot V$ in $\D$,  and 
$\O^h:V^h\ot V^h @>>> V^h\ot V^h$
is considered as a morphism in $\E$. 

\remark{Remark} We identify objects $V$, $V\ot \C[[h]]$
and $\C[[h]]\ot V$ in $\D$ as well as 
 $V^h$, $V^h\ot \X(\C[[h]])$
and $\X(\C[[h]])\ot V^h$ in $\E$.
\endremark

\medskip

\subhead  3. Classical Jacobi identity \endsubhead 
By the definition of Lie algebra, our distinguished morphism 
$T$ satisfies the antisymmetry property
$$
T\circ \s=-T\tag 3.1
$$
and the classical Jacobi identity
$$
T\circ (T\ot 1_V)=T \circ (1_V\ot T)\circ 
(1_{V\ot V\ot V}-\s\ot1_V)\tag 3.2
$$
The morphism $\s$ 
obviously satisfies the quantum Yang-Baxter equation.
$$
\s_{12}\circ\s_{23}\circ\s_{12}
=\s_{23}\circ\s_{12}\circ\s_{23}\tag 3.3
$$

\medskip

\subhead  4. Quantum Jacobi identity \endsubhead  Now let us 
apply the Kazhdan-Lusztig functor $\bold X$
to our distinguished morphisms $T$ and $\s$
and see what happens to the identities (3.1), (3.2) and (3.3).
Let us denote 
$$
M=(M_{V,V\ot V})\circ (1_{V^h}\ot M_{V,V}), 
\qquad
N=(M_{V\ot V,V})\circ (M_{V,V}\ot 1_{V^h})
$$
so $M,~~N:\X(V)\ot \X(V)\ot \X(V) @>>> \X(V\ot V\ot V)$ are
isomorphisms in $\E$. Note that by Proposition 25.8 of {\bf [6]}, 
we have
$$
\X(\Phi)=M\circ N^{-1}=(M-N)\circ N^{-1}+1 \tag 4.1
$$ 
where $\Phi$ here is considered as a morphism 
$V\ot V\ot V @>>>V\ot V\ot V$ in $\D$. We put
$$
\Psi=M^{-1}\circ N=M^{-1}\circ (N-M) +1\tag 4.2
$$ 
\proclaim {Theorem 4.1} The following identities hold for
morphisms in $\E$ 
$$
T^h\circ \s^h=-T^h,
\qquad
(\s^h)^2=1 
\tag 4.3 
$$
$$
T^h\circ (T^h\ot 1)=T^h \circ (1\ot T^h) 
\circ\Psi\circ(1_{V^h \ot V^h \ot V^h}
-\s^h\ot 1)\tag 4.4 
$$
$$
\Psi\circ\s^h_{12}\circ\Psi^{-1}
\circ\s^h_{23}\circ\Psi\circ\s^h_{12}=
\s^h_{23}\circ\Psi\circ\s^h_{12}
\circ\Psi^{-1}\circ\s^h_{23}\circ\Psi 
\tag 4.5
$$
\endproclaim
\demo{Proof} Follows from definitions, 
identities (3.1), (3.2) and (3.3),
and an elementary diagram search.\qed\enddemo

The identity (4.3) is our variant of quantum antisymmetry
and the identity (4.4) is our variant of 
quantum Jacobi identity.

\medskip

\subhead  5. The associator $\Psi$ \endsubhead 
It is well known (due to Le and Murakami 
{\bf [8]}; we use the exposition in {\bf [5]})
that the morphism $\Phi: V\ot V\ot V @>>> V\ot V\ot V$
can be presented as
$$
\Phi=1+\sum^{\infty}_{k=2}E_k(\O_{12},\O_{23})h^k\tag 5.1
$$
where $E_k(\O_{12},\O_{23})$ is a homogeneous polynomial
(of degree $k$) in $\O_{12}$ and $\O_{23}$. The explicit form
of the polynomials $E_k$ can be found in the references.
For example,
$$
E_2(\O_{12},\O_{23})=-\frac{\zeta (2)}{(2\pi\sqrt{-1})^2} 
[\O_{12},\O_{23}]=\frac{1}{24}[\O_{12},\O_{23}]\tag 5.2
$$
Let us now describe one formal polynomial map. Let $\C[A,B]$
be the ring of formal polynomials in two formal 
non-commuting variables over $\C$. 
We will construct a map 
$\alpha: \C[A,B] @>>> \C[C,D,Z,Z^{-1}]$, $\alpha(1)=1$
where $C$, $D$ and $Z$  are also formal pairwise 
non-commuting variables, $ZZ^{-1}=Z^{-1}Z=1$, as follows
$$
\alpha:A^{n_1}B^{m_1}A^{n_2}B^{m_2}\dots A^{n_j}B^{m_j}
\mapsto 
C^{n_1}Z^{-1}D^{m_1}Z \dots Z 
C^{n_j}Z^{-1}D^{m_j}Z 
$$
where $n_i$, $m_i\in \Bbb Z_{\geq 0}$, for $1\leq i\leq j$.
 We will now construct polynomials $G_k$ in variables
$\O^h_{12}$, $\O^h_{23}$, $\Psi$ and $\Psi^{-1}$ as follows:
consider $E_k(\O_{12},\O_{23})$ as a polynomial in
$A=\O_{12}$ and $B=\O_{23}$, take $\alpha(E_k)$,
and subsitute $\Psi$ and $\Psi^{-1}$ for $Z$ and $Z^{-1}$,
and $\O^h_{12}$ and  $\O^h_{23}$ for $C$ and $D$
(multiplication is composition of functors). Finally,
$F_k=-\Psi \circ G_k$.

\proclaim{Theorem 5.1}

a. 
$\Psi=1+\sum^{\infty}_{k=2}F_k(\O^h_{12},\O^h_{23},\Psi, 
\Psi^{-1})h^k$

b.
$\Psi^{-1}=
1+\sum^{\infty}_{k=2}G_k(\O^h_{12},\O^h_{23},\Psi, 
\Psi^{-1})h^k$.
\endproclaim

\demo{Proof} a. Notice that $\Psi^{-1}-1=N^{-1}\circ (M-N)$.
A diagram search shows that
$\X(E_k(\O_{12},\O_{23}))=
-M\circ F_k(\O^h_{12},\O^h_{23},\Psi, \Psi^{-1})\circ N^{-1}$.
This observation combined with (4.1), (4.2) and (5.1)
implies the claim by elementary calculations. 
Part b. is proved in an analogous way.
\qed\enddemo

\proclaim{Theorem 5.2} $\Psi, \Psi^{-1}\in 
\C[[h]][[\O^h_{12},\O^h_{23}]]$
\endproclaim

\demo{Proof} We will prove the claim by showing that
$\Psi \op{mod} h^n\in 
\C[h]/(h^n)[[\O^h_{12},\O^h_{23}]]$ and $\Psi^{-1}\op{mod} h^n\in 
\C[h]/(h^n)[[\O^h_{12},\O^h_{23}]]$  for any 
$n\in \Bbb Z_{>0}$. The latter is done by induction on $n$.
Indeed, $\Psi=1 \op{mod} h^2$, $\Psi^{-1}=1 \op{mod} h^2$, and if 
 $\Psi\op{mod} h^{n-2}$ and $\Psi^{-1}\op{mod} h^{n-2}$ 
are polynomials
in $\O^h_{12},\O^h_{23}$, then  
$$\matrix\format\l \\
\Psi
=1+\sum^{n-1}_{k=2}F_k(\O^h_{12},\O^h_{23},
\Psi \op{mod} h^{n-2}, \Psi^{-1} \op{mod} h^{n-2})h^k\hskip .1in
\op{mod} h^n \\
\Psi^{-1}=
1+\sum^{n-1}_{k=2}
G_k(\O^h_{12},\O^h_{23},\Psi\op{mod} h^{n-2}, 
\Psi^{-1}\op{mod} h^{n-2})h^k\hskip .1in\op{mod} h^n
\endmatrix
$$
are also polynomials in $\O^h_{12},\O^h_{23}$. Following 
this algorithm we can get an explicit
expression for  $\Psi\op{mod} h^n$ and  $\Psi^{-1}\op{mod} h^n$ 
for any 
$n\in \Bbb Z_{> 0}$.
\qed\enddemo

\medskip

\subhead 6. Acknowledgement \endsubhead I am very grateful 
to my advisor Igor Frenkel for the formulation of the problem 
and very useful discussions and to S. Majid, S. L. Woronowicz,
P. Etingof and V. Protsak for numerous useful discussions. 
I would like
to thank the organizers of the Program on Representation Theory
and Its Applications to Mathematical Physics at ESI 
for inviting me and giving me the opportunity to speak 
on the topic.
I am also grateful to ESI for their hospitality and 
financial support.

\enddocument
\end